\documentclass[12pt]{amsart}
\usepackage{amscd,times,fullpage}
\newtheorem{thm}{Theorem}
\newtheorem{prop}[thm]{Proposition}

\theoremstyle{remark}
\newtheorem{rem}[thm]{Remark}

\theoremstyle{definition}

\newtheorem{ex}[thm]{Example}

\newcommand{\C}{\mathbb{ C}}
\newcommand{\OO}{\mathbb{ O}}

\newcommand{\R}{\mathbb{ R}}
\newcommand{\HH}{\mathbb{ H}}

\newcommand{\Z}{\mathbb{ Z}}

\title{Geometric formality of homogeneous spaces and of biquotients}
\author{D.~Kotschick}
\address{Mathematisches Institut, Ludwig-Maximilians-Universit\"at
M\"unchen, Theresienstr.~39, 80333 M\"unchen, Germany}
\email{dieter@member.ams.org}
\author{S.~Terzi\'c}
\address{Faculty of Science, University of Montenegro, Dzordza Vasingtona bb, 81000 Podgorica, Montenegro}
\email{sterzic@cg.ac.yu}

\date{\today; \copyright{\ D.~Kotschick and Terzi\'c 2008}}


\begin{document}

\begin{abstract}
We provide examples of homogeneous spaces which are neither symmetric spaces nor real cohomology spheres, yet
have the property that every invariant metric is geometrically formal. We also extend the known obstructions to geometric
formality to some new classes of homogeneous spaces and of biquotients, and to certain sphere bundles.
\end{abstract}

\maketitle

\section{Introduction}\label{s:intro}

The notion of geometric formality was implicitly considered by Sullivan in the 1970s, see~\cite{DGMS,tokyo}, but
the systematic study of this concept began only comparatively recently~\cite{K}. A smooth manifold is geometrically formal if
it admits a Riemannian metric for which all exterior products of harmonic forms are harmonic. Such a metric is
then also called formal. Geometric formality clearly implies formality in the
sense of Sullivan, and is even more restrictive. As compact
symmetric spaces are the classical examples of geometrically formal
manifolds, it is natural to explore this notion in the context of
homogeneous spaces, or, more generally, of manifolds with large symmetry groups. 

Trying to come up with generalizations of symmetric spaces, one might think first 
of isotropy irreducible spaces. These are the homogeneous spaces $G/H$ for which the
isotropy representation of $H$ on $T_{eH}(G/H)$ is irreducible. Such a space is
strongly isotropy irreducible if the restriction of the isotropy
representation to the identity component of $H$ is also irreducible.
These manifolds were originally classified by Manturov, and were
further studied by Wolf and others, cf.~\cite{Besse}. They share many
properties of symmetric spaces, and indeed irreducible symmetric
spaces are isotropy irreducible. A conceptual relationship between
symmetric spaces and isotropy irreducible ones is explained in~\cite{WZ}.
However, the similarities between symmetric spaces and isotropy irreducible
ones do not extend to (geometric) formality. Indeed, there are a
number of strongly isotropy irreducible spaces which, by the results
of~\cite{Doan Kuin'} are not of Cartan type, and, therefore~\cite{KT}, are not
formal in the sense of Sullivan. {\it A fortiori}, they cannot be
geometrically formal.
\begin{ex}
The compact homogeneous spaces $SU(pq)/(SU(p)\times SU(q))$ for
$p,q\geq 3$, $SO(78)/E_{6}$, and $SO(n^{2}-1)/SU(n)$ for $n\geq 3$ are
strongly isotropy irreducible, but are not formal in the sense of
Sullivan.
\end{ex}

Another class of homogeneous spaces generalizing the symmetric ones 
consists of the so-called generalized symmetric spaces, sometimes called $k$-symmetric.
These are defined by replacing the involution in the definition of symmetric 
spaces by a symmetry of order $k$, see~\cite{W-G,T}.
In~\cite{KT} we proved that all generalized symmetric spaces of
compact simple Lie groups are formal in the sense of Sullivan, and
that many of them are \emph{not} geometrically formal.

The main purpose of this paper is to prove that, in spite of all these negative results,
there are indeed homogeneous spaces which are neither (homotopy equivalent to)
symmetric spaces nor products of real homology spheres (which are trivially geometrically formal~\cite{K}),
yet are geometrically formal. We shall prove the following:
\begin{thm}\label{t:main}
All homogeneous metrics on the following homogeneous spaces are geometrically formal:
\begin{enumerate}
\item the real Stiefel manifolds $V_{4}(\R^{2n+1}) = SO(2n+1)/SO(2n-3)$ for $n\geq 3$,
\item the real Stiefel manifolds $V_3(\R^{2n})=SO(2n)/SO(2n-3)$ for $n\geq 3$,
\item the complex Stiefel manifolds $V_2(\C^n) =SU(n)/SU(n-2)$, for $n\geq 5$,
\item quaternionic Stiefel manifolds $V_2(\HH^n)=Sp(n)/Sp(n-2)$, for $n\geq 3$,
\item the octonian Stiefel manifold $V_2(\OO^2)=Spin(9)/G_2$, and
\item the space $Spin(10)/Spin(7)$.
\end{enumerate}
Moreover, none of these spaces is homotopy equivalent to a symmetric space.
They are not homotopy equivalent to products of real cohomology spheres, except possibly for $V_3(\R^{2n})$ with $n$ even.
\end{thm}
The space $Spin(10)/Spin(7)$ in the theorem corresponds to a non-standard embedding of $Spin(7)$ in $Spin(10)$, so that the 
quotient is not $V_3(\R^{10})$, but a manifold with the same real and different integral homology. The homology of this space,
unlike that of $V_3(\R^{10})$, is torsion-free. 

In Section~\ref{s:new} we exhibit a very simple mechanism to prove that $G$-invariant metrics on certain homogeneous spaces
$G/H$ with simple cohomology rings are geometrically formal. This mechanism in fact gives a new proof that certain symmetric spaces
are geometrically formal, without using the symmetric space structure, but only the description of the cohomology ring. 
The argument applies to all homogeneous spaces which have the real cohomology of a product of odd-dimensional spheres.
The examples (3)-(6) in Theorem~\ref{t:main} are all the homogeneous spaces which have the integral cohomology of such a 
product, but are not obviously diffeomorphic to products or to symmetric spaces. There are many more examples with the real
cohomology of such a product, but different integral cohomology. We discuss the two infinite sequences of real Stiefel manifolds occurring in (1) 
and (2), leaving aside the other, sporadic, examples.
Also in Section~\ref{s:new}, by considering the Alof-Wallach spaces~\cite{AW},
we show that our results do not extend to all homogeneous spaces with the cohomology
algebra of a product of spheres, if one does not insist that the spheres be odd-dimensional. 
In Section~\ref{s:at} we show that the homogeneous spaces listed in
Theorem~\ref{t:main} are not homotopy equivalent to symmetric spaces or to non-trivial products, thereby completing
the proof of the theorem. 

In the final two sections of this paper we add to the negative results of~\cite{K,KT,KT2} by providing 
further examples of manifolds which, though formal in the sense of Sullivan, are not geometrically formal. All the examples 
we give here are simply connected and of dimension six. In Section~\ref{s:biq} we consider certain classes of biquotients 
in the sense of Eschenburg~\cite{E0,E}, and in Section~\ref{S2-bundles} we consider two-sphere bundles over $\C P^2$. 
Many of these two-sphere bundles are known to carry special metrics of cohomogeneity one by the results of Grove and 
Ziller~\cite{GZ}. Both these collections of examples generalize
the discussion of the flag manifold $SU(3)/T^{2}$ carried out in~\cite{KT}. A different
generalization, to certain partial flag manifolds of higher dimensions, is contained in~\cite{KT2}.

\subsection*{Acknowledgement}
The work of the first author was supported in part by the DFG Priority Program in Global Differential Geometry 
and by The Bell Companies Fellowship at the Institute for Advanced Study in Princeton.

\section{Examples of geometrically formal homogeneous spaces}\label{s:new}

In this section we describe a class of homogeneous spaces for which
any homogeneous metric is formal. 

First we recall the well known fact that on a compact homogeneous 
space $G/H$ the harmonic forms of any homogeneous metric are invariant.
To check this, let $h$ be a harmonic form with respect to some homogeneous metric on $G/H$.
Then $h$ and $\ast h$ can be written as $h = h_{i} + d\alpha$,
$\ast h = (\ast h)_{i} +d\beta$, where $h_{i}=\int _{G}g^{*}h$ and
$(\ast h)_{i}=\int _{G}g^{*}(\ast h)$ with respect to Haar measure
with total volume $1$. Since $h_{i}$ and $(\ast h)_{i}$
are invariant forms we have
$$
  \ast h_{i} = \ast \int _{G}g^{*}h = \int _{G}\ast (g^{*}h)
  = \int _{G}g^{*}(\ast h) = (\ast h)_{i} \ .
$$
It follows that $d(\ast h_{i}) = d(\ast h)_{i} = 0$ and, thus,
$h_{i}$ is harmonic. This means that $h = h_{i}$, i.e. $h$ is
an invariant form.

As a first application of this fact we have the following:
\begin{prop}\label{p:1}
Let $G$ be a compact connected Lie group and $H$ a closed connected
subgroup with the property that $G/H$  is of even dimension $2k$, and all the 
real cohomology is in degrees $0$, $k$ and $2k$. Then any homogeneous metric on $G/H$ is formal.
\end{prop}
\begin{proof}
Because of the cohomology structure of $G/H$, to prove that $G/H$ is
geometrically formal with respect to some metric, it is enough to prove that
$x\wedge y$ is a harmonic form, for any two harmonic $k$-forms $x$ and $y$.
If the metric $g$ is homogeneous, then
according to the previous observation, $x\wedge y$ is an invariant form and being of top
degree we have that $x\wedge y=c \cdot dvol$, where $c$ is constant. Thus
$x\wedge y$ is harmonic and $g$ is formal.
\end{proof}
\begin{ex}\label{ex:11}
This Proposition applies in particular to the following spaces: the complex projective plane
$\C P^2 = SU(3)/S(U(2)\times U(1))$, the quaternionic projective plane $\HH P^2 = Sp(3)/(Sp(2)\times Sp(1))$,
the Cayley plane $\OO P^2 = F_4/Spin(9)$, and $G_2/SO(4)$, cf.~\cite{BH}. 
\end{ex}
These spaces are all symmetric, but the argument proving geometric formality does not use the 
symmetric space structure. Unfortunately there are no non-symmetric homogeneous spaces 
to which we could apply Proposition~\ref{p:1}. Indeed, if the number $k$ in Proposition~\ref{p:1} is even, 
then such spaces belong to the class of rank one homogeneous spaces in the terminology of~\cite{OR}, 
while for odd $k$ they belong to the class of rank two homogeneous spaces. Examining the classification 
of homogeneous spaces of rank one and two given by Onishchik~\cite{OR}, one sees that there are no
examples other than the symmetric spaces mentioned in Example~\ref{ex:11}.  

To get new examples we need the following slight variation of Proposition~\ref{p:1}:
\begin{prop}\label{prod}
Let $G$ be a compact connected Lie group and $H$ a closed connected
subgroup with the property that $H^{*}(G/H;\R)=\wedge (x, y)$, where $x$ and $y$
are of odd degrees. Then any homogeneous metric on $G/H$ is formal.
\end{prop}
The spaces listed in Theorem~\ref{t:main} all have cohomology rings of this form, cf.~\cite{On,OR,Kramer}.
Thus Proposition~\ref{prod} shows that those spaces are geometrically formal.
\begin{proof}
Because of the cohomology structure of $G/H$, to prove that $G/H$ is
geometrically formal with respect to some metric, it is enough to prove that
$x\wedge y$ is a harmonic form, where $x$ and $y$ are the harmonic representatives
of the cohomology generators. If the metric $g$ is homogeneous, then
as before, $x\wedge y$ is an invariant form and being of top
degree we have that $x\wedge y=c \cdot dvol$, where $c$ is constant. 
It follows that $x\wedge y$ is harmonic and $g$ is formal.
\end{proof}

\begin{rem}\label{nonprod}
At this point it is useful to recall that if $G$ is simple and we endow
$G/H$ with the submersion metric $g$ of a biinvariant metric on $G$,
then the Riemannian homogeneous space $(G/H, g)$ is irreducible as a
Riemannian manifold~\cite{KN}. This means that such a $g$ cannot be a product
metric. These are examples of normal homogeneous metrics.
\end{rem}
As a simple application of our discussion so far, we can prove the following.
\begin{prop}\label{st-prod}
The group $SU(4)$ acts transitively on $S^{5}\times S^{7}$. All $SU(4)$-homogeneous metrics
for this action are formal. Furthermore, the normal homogeneous metrics are not symmetric.
\end{prop}
\begin{proof}
Firstly, it is clear that $SU(4)$ acts 
transitively on $S^7$, with isotropy group $SU(3)$. Secondly, $SU(4)$ acts transitively on $S^5$
via the double covering $SU(4)\longrightarrow SO(6)$. The isotropy group of this action is the preimage
of $SO(5)$ under the covering, which can be identified with $Sp(2)$. Now take the product action
on $S^{5}\times S^{7}$. This is still transitive, for example because the restriction of the action of 
$SU(4)$ on $S^5$ to $SU(3)$ is still transitive. It follows that $S^{5}\times S^{7}$ is a homogeneous
space of $SU(4)$ with isotropy group $Sp(2)\cap SU(3) = SU(2)$.
 
Proposition~\ref{prod} implies that all $SU(4)$-homogeneous metrics on $S^5\times S^7$ are formal.
By Remark~\ref{nonprod}, it follows that we cannot get the normal homogeneous metrics as product metrics. 
In particular the metric on $S^5\times S^7$ which is the product of the symmetric space metrics on the factors, 
though also formal, is not normal homogeneous for the $SU(4)$-action.
\end{proof}

To end this section, we want to show that Proposition~\ref{prod} is 
sharp in the sense that it does not extend to arbitrary homogeneous spaces with a 
cohomology algebra of the form $\wedge (x, y)$ with $x$ of even degree and $y$ of odd degree.

A convenient class of examples to consider for this purpose are the so-called 
Aloff-Wallach spaces. These are homogeneous spaces of the form $N_{k,l}=SU(3)/T^{1}$,
where $T^{1}$ is embedded as the diagonal matrices $D(z^{k},z^{l},z^{-k-l})$
with $k$ and $l$ coprime integers with $kl(k+l)\neq 0$. Obviously, they are all
homogeneous spaces of Cartan type, and are therefore~\cite{KT} formal in the sense of Sullivan. It is
also easy to see that they all have the real cohomology of $S^{2}\times S^{5}$.
The name derives from~\cite{AW}, where Aloff and Wallach proved that
these spaces have homogeneous metrics of positive sectional
curvature. The numerical conditions on $k$ and $l$ are there to make
sure that $T^{1}$ acts on $\C^{3}$ without nonzero fixed points.

Since $SU(3)$ endowed with a biinvariant metric is geometrically formal, it is
natural to consider the submersion metrics on $SU(3)/T^{1}$ which we get from
the principal circle fibration
\begin{equation}\label{pr}
T^{1}\longrightarrow SU(3)\stackrel{\pi}{\longrightarrow} N_{k,l} \ .
\end{equation}
Such a metric is often called a normal homogeneous metric.
\begin{thm}
The normal homogeneous metrics on Aloff-Wallach spaces are not formal.
\end{thm}
\begin{proof}
Assume that $M=SU(3)/T^{1}$ is geometrically formal for some embedding
$T^{1}\subset SU(3)$. Since $M$ has the real cohomology of $S^{2}\times S^{5}$,
it carries two harmonic forms $\omega_{2}$ and $\omega_{5}$ such that
$\omega _{2}^{2}=0$ and $\omega_{2}\wedge \omega_{5}$ is a volume form.
Since $SU(3)$ is simply connected, the Euler class $e$ of the principal
bundle~\eqref{pr} is not zero, so that $e=\lambda [\omega _{2}]$ with
$\lambda \neq 0$. We now normalise our form so that $\lambda =1$. In
the following calculation we also ignore non-zero constants.

There exists a connection form $\alpha$ on the principal bundle~\eqref{pr} such that
\begin{equation}\label{con}
\pi^{*}(\omega_{2})=d\alpha \ ,
\end{equation}
where $\pi$ is the projection in~\eqref{pr}.

Let $\eta_{5}=\pi^{*}(\omega_{5})$ and $\eta_{2}=\pi^{*}(\omega_{2})$.
Then $\eta _{2}$ and $\eta _{5}$ are closed and
\begin{equation}\label{sq}
      \eta _{2}^{2}=0 \ .
      \end{equation}
Also,
\begin{equation}
      \ast \eta _{5}=\ast \pi ^{*}(\omega _{5})=
      \alpha \wedge \pi ^{*}(\ast \omega _{5})=\alpha \wedge \eta _{2} \ .
      \end{equation}
Then~\eqref{con} gives
\begin{equation}
      d(\ast \eta _{5})=d\alpha \wedge \eta _{2}=\eta _{2}^{2}=0 \ .
      \end{equation}
This implies that $\eta _{3}=\ast \eta _{5}$ is a harmonic form on $SU(3)$
with the biinvariant metric. Since the harmonic forms on $SU(3)$ coincide
with the biinvariant ones, we get that $\eta _{3}$ has the form
\[
\eta _{3}(X,Y,Z)=\left< X,[Y,Z]\right> \; \mbox{for} \; X,Y,Z\in
\mathfrak{su}(3) \ .
\]
On the other hand we have a natural direct sum decomposition
\begin{equation}\label{dc}
      \mathfrak{su}(3)=\mathfrak{t}^{1}\oplus
      (\mathfrak{su}(3)/\mathfrak{t}^{1}) \ ,
      \end{equation}
and~\eqref{sq} implies that there exists $5$-dimensional subspace
$\mathcal{K}$ in $\mathfrak{su}(3)/\mathfrak{t}^{1}$, such that
$i_{v}(\eta _{3})=0$ for any vector $v\in \mathcal{K}$.

Let $H_{1},H_{2},E_{1},E_{2},F_{1},F_{2}$ be canonical (Chevalley) generators
of the Lie algebra $\mathfrak{su}(3)$. Consider the subspace
$\mathcal{L}$ spanned by the vectors $H_{1}, H_{2}, E_{1} , F_{1}$. We are going
to show that for any $X\in \mathcal{L}$, $i_{X}(\eta _{3})\neq 0$.
Any $X\in \mathcal{L}$ can be written in the form
$X = aH_{1} + bH_{2} + cE_{1} + dF_{1}$. To prove the above we consider the
following cases.

{\bf 1.} If $d\neq 0$, then using the well-known relation between canonical generators
of a simple Lie algebra and the root spaces related to the Killing
form, we get $\eta _{3}(E_1, H_1, X) = -2d\langle E_1, F_1\rangle \neq 0$.

{\bf 2.} If $d = 0$ and $c\neq 0$, then we have $\eta _{3}(F_1, H_1, X) = 2c\langle F_1, E_1\rangle \neq 0$.

{\bf 3.} If $c = d = 0$, then we have
$\eta _{3}(F_1, X, E_1) = (2a-b)\langle F_{1}, E_{1}\rangle$ and
$\eta _{3}(F_2, X, E_2) = (2b-a)\langle F_{2}, E_{2}\rangle$, so there always are $E_i, F_i$ for which $\eta _{3}(F_i, X, E_i)\neq 0$.

The above implies that $\mathcal{L}\cap \mathcal{K} = 0$, which is impossible for dimension reasons.
\end{proof}

\begin{rem}
There are exactly two fibrations with fiber $S^{2}$ over $S^{5}$,
cf.~section~\ref{S2-bundles} below. The trivial bundle $S^{2}\times S^{5}$ is
of course geometrically formal with respect to product metrics. The
nontrivial bundle is the $3$-symmetric space $SU(3)/T^{1}$, where
$T^{1}$ is embedded inside an $SU(2)\subset SU(3)$. In the above
notation, this is the case $k=-l$, which is excluded in the
definition of Aloff-Wallach spaces. Nevertheless, the above argument
applies to show that a normal homogeneous metric is not geometrically
formal.
\end{rem}

\begin{rem}
The Aloff-Wallach spaces have interesting homogeneous metrics which are
not normal homogeneous. These include metrics of positive sectional
curvature~\cite{AW} and Einstein metrics, some of which admit Killing spinors. The
latter metrics are not geometrically formal because of the following result,
communicated to us by U.~Semmelmann~\cite{Sem}:
    A metrically formal Riemannian spin manifold $M$ of dimension
    $\geq 5$ admitting a nontrivial Killing spinor must have vanishing second Betti number.
    
It is known that a metric admitting a nontrivial Killing spinor must
be Einstein, and thus is very special.
\end{rem}

\section{Some algebraic topology of Stiefel manifolds}\label{s:at}

In this section we prove that the homogeneous spaces listed in Theorem~\ref{t:main} are not homotopy equivalent to symmetric spaces or to 
products of real homology spheres, except for the second property in the case of $V_3(\R^{2n})$ with $n$ even. Together with Proposition~\ref{prod}
proved in the previous section, this completes the proof of Theorem~\ref{t:main}.

Let us consider first the complex Stiefel manifolds $V_2(\C^n)=SU(n)/SU(n-2)$ consisting of orthonormal pairs of vectors for the 
standard Hermitian inner product on $\C^n$, $n\geq 3$. Projecting such a pair to its first entry we obtain a smooth fibration of $V_2(\C^n)$
over $S^{2n-1}$ with fibre $S^{2n-3}$. It is a classical problem in homotopy theory to determine when the total space of such a 
fibre bundle is homotopy equivalent to the product of base and fibre. 
In the case at hand, for even $n$ the action of the quaternions on $\C^n = \HH^{n/2}$
defines a section of the fibration, which splits the long exact homotopy sequence and makes $V_2(\C^n)$ indistinguishable 
from $S^{2n-1}\times S^{2n-3}$ at the level of homotopy groups. 
However, the following result was proved by James and Whitehead~\cite{JW}
modulo Adams's solution of the Hopf invariant one problem:
\begin{thm}{\rm (James--Whitehead~\cite{JW}, Adams~\cite{Adams})}\label{t:JW}
For $n\geq 3$, the Stiefel manifold $V_2(\C^n)=SU(n)/SU(n-2)$ is not homotopy equivalent to $S^{2n-1}\times S^{2n-3}$ unless possibly if $n=4$.
\end{thm}
In fact, James and Whitehead~\cite{JW} proved that if $V_2(\C^n)=SU(n)/SU(n-2)$ is homotopy equivalent to $S^{2n-1}\times S^{2n-3}$,
then $\pi_{4n-1}(S^{2n})$ contains an element of Hopf invariant one. By Adams's result~\cite{Adams} it follows that $n\in\{1,2,4\}$. 
This combination of the results of~\cite{JW} and~\cite{Adams} is also mentioned in~\cite[Theorem 1.7]{James}, 
but the exceptional case is misstated there. In the notation of~\cite{James}, the exceptional case should be denoted $n=4$, $k=2$, and not $n=k=2$, 
which is what is written in~\cite{James}. 

\begin{rem}
Some years after the results of James--Whitehead~\cite{JW} and Adams~\cite{Adams}, Gilmore~\cite{G} showed that $\pi_{2n-2}(V_2(\C^n))$ is 
trivial for odd $n$. As $\pi_{2n-2}(S^{2n-1}\times S^{2n-3})=\pi_{2n-2}(S^{2n-3})=\Z_2$, this gives another proof of Theorem~\ref{t:JW} for 
odd $n$. As we remarked earlier, there is no such proof in the case of even $n$.
\end{rem}

The arguments of James and Whitehead only show that the existence of an element of Hopf 
invariant one is necessary for $V_2(\C^n)=SU(n)/SU(n-2)$ to be homotopy equivalent to $S^{2n-1}\times S^{2n-3}$. It turns out 
that this condition is in fact sufficient not just for homotopy equivalence, but for diffeomorphism.
This following result completes the proof of case (3) in Theorem~\ref{t:main}.
\begin{thm}\label{t:cx}
If a complex Stiefel manifold $V_2(\C^n)=SU(n)/SU(n-2)$ with $n\geq 3$ is homotopy equivalent to a symmetric space or to a product of real 
homology spheres, then $n=3$ or $4$. In the first case $V_2(\C^3)=SU(3)$ is a symmetric space not homotopy equivalent to a product of 
real homology spheres. In the second case $V_2(\C^4)=SU(4)/SU(2)$ is diffeomorphic to $S^5\times S^7$.
\end{thm}
\begin{proof}
Suppose that $V_2(\C^n)$ is homotopy equivalent to a product $X_1\times X_2$ of real homology spheres. Then because $V_2(\C^n)$ is simply
connected, so are both $X_i$. Moreover, because $V_2(\C^n)$ has the integral homology of a product of spheres, it follows that each $X_i$
is an integral homology sphere. Thus the $X_i$ are homotopy spheres. 

Now, if $V_2(\C^n)$ is homotopy equivalent to a product of (homotopy) spheres, then Theorem~\ref{t:JW} gives $n=4$. Conversely, $V_2(\C^4)$ is an $S^5$-bundle 
over $S^7$ with structure group $SU(4)$. The corresponding clutching map gives an element of $\pi_6(SU(4))$, which is trivial~\cite{BS}. Thus 
$V_2(\C^4)$ is diffeomorphic to the total space of the trivial bundle $S^5\times S^7$, compare Proposition~\ref{st-prod}.

Next suppose that $V_2(\C^n)$ is homotopy equivalent to a symmetric space. Onishchik~\cite{On,OR} and Kramer~\cite{Kramer} have classified 
the homogeneous spaces with the cohomology of a product of odd-dimensional spheres. If we assume that $G/H$ has the integral cohomology of 
such a product, then, apart from products of spheres and a handful of symmetric spaces of dimension $\leq 14$, the remaining examples are 
exactly the spaces listed as (3)--(6) in Theorem~\ref{t:main} and the symmetric space $E_6/F_4$ with the cohomology of $S^9\times S^{17}$.
If $n\geq 5$, then the dimension of $V_2(\C^n)$ is $\geq 16$, 
so that we do not have to consider the sporadic irreducible symmetric spaces of dimension $\leq 14$ arising in the classification of 
Onishchik~\cite{On,OR} and Kramer~\cite{Kramer}. Similarly, $V_2(\C^n)$ can not be homotopy equivalent to $E_6/F_4$ for dimension reasons.
If $V_2(\C^n)$ were homotopy equivalent to a reducible symmetric space, then each factor 
would be a homotopy sphere, and we would be in the case excluded already.

Finally, $V_2(\C^3)=SU(3)$ is of course a symmetric space, but it is not homotopy equivalent to a product of spheres. For if it were, then
the spheres would have dimensions $3$ and $5$, leading to a contradiction with $\pi_6(SU(3))=\Z_6$, cf.~\cite{BS}, since $\pi_6(S^3)=\Z_{12}$.
\end{proof}

In a similar way one proves the corresponding result for the quaternionic Stiefel manifolds. The statement is cleaner, as an exceptional
case not excluded by the Hopf invariant does in fact not occur. This following result completes the proof of case (4) in Theorem~\ref{t:main}.
\begin{thm}
No quaternionic Stiefel manifold $V_2(\HH^n)=Sp(n)/Sp(n-2)$ with $n\geq 2$ is homotopy equivalent to a product of real 
homology spheres. If it is homotopy equivalent to a symmetric space, then $n=2$.
\end{thm}
\begin{proof}
If $V_2(\HH^n)$ is homotopy equivalent to a non-trivial product, then, as in the previous proof, we may assume that each factor is a sphere. 
Assuming that $V_2(\HH^n)$ is homotopy equivalent to $S^{4n-5}\times S^{4n-1}$, James and Whitehead~\cite[Theorem 1.21]{JW} proved 
that $\pi_{8n-1}(S^{4n})$ contains an element of Hopf invariant one. The result of Adams~\cite{Adams} then implies $n=2$. However,
 for $n=2$, we have $V_2(\HH^2)=Sp(2)$, and this is not homotopy equivalent to a product of spheres, for example because 
 $\pi_6(Sp(2))$ is trivial as first proved by Borel and Serre~\cite{BS}.

Next suppose that $V_2(\HH^n)$ is homotopy equivalent to a symmetric space. If $n\geq 3$, then the dimension of $V_2(\HH^n)$ is 
$\geq 18$, so that, again, we do not have to consider the sporadic irreducible symmetric spaces of dimension $\leq 14$ arising in the classification due to
Onishchik~\cite{On,OR} and Kramer~\cite{Kramer}. For dimension reasons $E_6/F_4$ can not occur either.
In the reducible case $V_2(\HH^n)$ would be homotopy equivalent to a product of spheres, and this we have excluded already.
Finally, for $n=2$ we have the symmetric space $V_2(\HH^2)=Sp(2)$.
\end{proof}
\begin{rem}
Instead of using the results of James--Whitehead~\cite{JW} and Adams~\cite{Adams}, we could, in most cases, appeal to the calculations of Oguchi~\cite{O}.
On the one hand we have $\pi _{4n-2}(S^{4n-1} \times S^{4n-5})=\pi _{4n-2}(S^{4n-5})=\Z_{24}$ for $n\geq 3$. On the other, by~\cite{O}, $\pi_{4n-2}(V_2(\HH^n))=\Z_d$,
where $d=gcd\{n,24\}$. Thus, whenever $n$ is not divisible by $24$, one concludes that $V_2(\HH^n)$ is not homotopy equivalent to a product of spheres.
\end{rem}

To complete the proof of case (5) in Theorem~\ref{t:main}, we prove the following:
\begin{thm}
The octonian Stiefel manifold $V_2(\OO^2)=Spin(9)/G_2$ is not homotopy equivalent to a symmetric space, or to a product of real homology spheres.
\end{thm}
\begin{proof}
This manifold has the integral homology of $S^7\times S^{15}$. As in the previous proofs, if $V_2(\OO^2)$ were homotopy equivalent
to a non-trivial product, then this product would have to be $S^7\times S^{15}$. We now prove that $V_2(\OO^2)$ and $S^7\times S^{15}$
are in fact not homotopy equivalent, distinguished by their homotopy groups.

On the one hand we have
$$
\pi_{14}(S^7\times S^{15}) = \pi_{14}(S^7) = \Z_{120} \ .
$$
On the other hand we will show that $\pi_{14}(V_2(\OO^2))$ has order at most $16$. For this we consider the principal $G_2$-bundle 
with total space $Spin(9)$ over $V_2(\OO^2)$, and the following piece of its exact homotopy  sequence:
$$
\ldots\longrightarrow\pi_{14}(Spin(9))\longrightarrow \pi_{14}(V_2(\OO^2)) \longrightarrow \pi_{13}(G_2) \longrightarrow\ldots
$$
The group on the right is trivial and the group on the left is $\Z_8\oplus\Z_2$, by calculations of Mimura~\cite{M}.

The proof that $V_2(\OO^2)$ is not homotopy equivalent to a symmetric space is the same as for the complex and quaternionic Stiefel manifolds.
First of all, reducible symmetric spaces cannot arise because $V_2(\OO^2)$ is not homotopy equivalent to a product, as we just proved. Second of all,
there is no irreducible symmetric space with the correct cohomology.
\end{proof}

Next we deal with the space $Spin(10)/Spin(7)$, case (6) in Theorem~\ref{t:main}.
\begin{thm}\label{t:ko}
The quotient $Spin(10)/Spin(7)$ is not homotopy equivalent to a symmetric space, or to a product of real homology spheres.
\end{thm}
\begin{proof}
Consider the principal $Spin(9)$-bundle $Spin(10)\longrightarrow Spin(10)/Spin(9)=S^9$. The spin representation of $Spin(9)$ on $\R^{16}$
associates to this principal bundle a real vector bundle $V$ of rank $16$. Our homogeneous space $X=Spin(10)/Spin(7)$ is the unit sphere bundle
in $V$, with fibre $Spin(9)/Spin(7)=S^{15}$. 

The manifold $X$ has the integral homology of $S^9\times S^{15}$. As in the previous proofs, if $X$ were homotopy equivalent
to a non-trivial product, then this product would have to be $Y=S^9\times S^{15}$. We will show below that in fact $X$ and $Y$
are not homotopy equivalent. 
The proof that $X$ is not homotopy equivalent to a symmetric space is then the same as for the complex and quaternionic Stiefel manifolds.
First of all, reducible symmetric spaces cannot arise because $X$ is not homotopy equivalent to a product. Second of all,
there is no irreducible symmetric space with the correct cohomology.

Now we begin the proof that $X$ and $Y=S^9\times S^{15}$ are not homotopy equivalent.
The principal bundle $Spin(10)\longrightarrow S^9$ corresponds to an element of $\pi_8(Spin(9))$,
whose image in $\pi_8(SO(16))$ classifies $V$. As $V$ is non-trivial, this element is the non-trivial element of $\pi_8(SO(16))=\Z_2$.

Recall that the fibration $X\longrightarrow S^9$ has a section for dimension reasons. Now a result of James and 
Whitehead~\cite[Cor.~(1.9)]{JW} tells us that $X$ is homotopy equivalent to $Y=S^9\times S^{15}$ if and only if a certain invariant 
$\lambda(X)$ vanishes. This invariant is an element in the group $\Lambda_{9,15}=Im J/Im P$, where 
$$
J \colon \pi_8(SO(15))\longrightarrow \pi_{23}(S^{15})
$$
is the classical $J$-homomorphism and 
$$
P\colon\pi_9(S^{15})\longrightarrow \pi_{23}(S^{15})
$$
is, in our case, the zero-map, as its domain is zero. Thus, in our case, $\lambda(X)\in J( \pi_8(SO(15))$, and we only have to determine whether 
this is zero, or not. Unravelling the definition of $\lambda$ given in~\cite{JW}, we find the following: if we identify $\pi_8(SO(15))$
with $\pi_8(SO(16))$, then $\lambda(X)$ is just the image under the classical $J$-homomorphism of the classifying element of $V$ in $\pi_8(SO(16))$.
Now $\pi_8(SO(15))$ has order $2$, and we know that the image of the $J$-homomorphism in this degree is the group $J(S^9)$, also of order $2$,
cf.~\cite{Bott,Huse}. Thus the $J$-homomorphism is an isomorphism, and $\lambda(X)\neq 0$ because the classifying map of $V$
represents the non-zero element of $\pi_8(SO(16))$.
\end{proof}

\begin{rem}
The existence of a section to the fibration $S^{15}\longrightarrow X\stackrel{\pi}{\longrightarrow} S^9$ implies that $X$ and $Y=S^9\times S^{15}$ 
have isomorphic homotopy groups. Thus there can be no easy argument to prove that they are not homotopy equivalent.
\end{rem}


It remains to discuss the homotopy types of the real Stiefel manifolds occurring in cases (1) and (2) of Theorem~\ref{t:main}.
They are geometrically formal by Proposition~\ref{prod}, and they are not homotopy equivalent to 
products of spheres because of the presence of torsion in their integral homology. However, for these manifolds it is more difficult to 
exclude homotopy equivalence to products of real homology spheres, because the factors in such a decomposition would not necessarily be 
homotopy spheres. The following result completes the proof of Theorem~\ref{t:main}.
\begin{thm}
The real Stiefel manifolds $V_4(\R^{2n+1})$ and $V_3(\R^{2n})$ with $n\geq 3$ are not homotopy equivalent to symmetric spaces.
They are not homotopy equivalent to non-trivial products, except possibly for $V_3(\R^{2n})$ with $n$ even.
\end{thm}
\begin{proof}
These manifolds have the real cohomology algebras of $S^{4n-5}\times S^{4n-1}$ and $S^{2n-1}\times S^{4n-5}$ respectively.
Therefore, if one of them is homotopy equivalent to a non-trivial product, then the factors are real homology spheres of dimensions 
$4n-5$ and $4n-1$, respectively $2n-1$ and $4n-5$.

The cohomology with coefficients in $\Z_2$, and the Steenrod operations on it, for the real Stiefel manifolds was determined by
Borel~\cite{Borel3}. From the structure of the Steenrod operations, Hsiang and Su~\cite{HS} deduced that $V_4(\R^{2n+1})$ and $V_3(\R^{2n})$ can not be homotopy 
equivalent to non-trivial products, except possibly in the second case when $n$ is even.

Suppose now that one of these manifolds is homotopy equivalent to a symmetric space. Then this symmetric space is reducible, and we 
are in the exceptional case where $V_3(\R^{2n})$ could be homotopy equivalent to a product. To see this, consider the irreducible 
symmetric spaces. The computations of their real cohomology algebras~\cite{Borel,Takeuchi}, or of their rational homotopy 
groups~\cite{T1}, show that among them only $SU(5)/SO(5)$, $SU(6)/Sp(3)$, $E_{6}/F_{4}$, $Sp(2)$, $SU(3)$, $G_2$ and $Spin(4)$ have 
the real cohomology algebra of a product of odd-dimensional spheres, with the dimensions given by
the pairs $(5, 9)$, $(5, 9)$, $(9,17)$, $(3,7)$, $(3,5)$, $(3,11)$ and $(3,3)$ respectively. None of these pairs of dimensions is of the form 
$(4n-5, 4n-1)$ or $(2n-1, 4n-5)$ with $n\geq 3$. Therefore, $V_4(\R^{2n+1})$ or $V_3(\R^{2n})$ can not be homotopy 
equivalent to any irreducible compact simply connected symmetric space.

Finally we have to consider the possibility that $V_3(\R^{2n})$ with $n$ even could be homotopy equivalent to a reducible symmetric 
space. Then each factor would be a compact symmetric space with the real cohomology of $S^{2n-1}$, respectively $S^{4n-5}$. 
As $V_3(\R^{2n})$ has torsion in its integral homology at least one of the factors has to be a simply connected symmetric space which has the real cohomology of a 
sphere, but different integral cohomology. According to the classification, due to Onishchik~\cite{On,OR}, of homogeneous spaces whose real cohomology is that of 
a sphere, the only possibility is $SU(3)/SO(3)$, of dimension $5$. For even $n$, neither $2n-1$ nor $4n-1$ can equal $5$, and this contradiction completes the proof.
\end{proof}

\begin{rem}
Recall that projection of a $3$-frame onto its first entry defines a smooth fibration of $V_3(\R^{2n})$ over $S^{2n-1}$, with fibre $V_2(\R^{2n-1})$.
Note that the fibre $V_2(\R^{2n-1})$ is a real homology sphere, but not an integral one. 
As far as we know, it is still an open problem whether for some even $n$ the total space could be homotopy equivalent, or even diffeomorphic, to 
the product of base and fibre, cf.~the discussion in~\cite{James}. If this happens for some $n$, then by a result of James~\cite{James}, this $n$ is not just even, but a 
power of $2$.
\end{rem}

Onishchik~\cite{On,OR} and Kramer~\cite{Kramer} have classified all the compact homogeneous spaces with the real cohomology of a product of 
odd-dimensional spheres. In addition to non-trivial products and the spaces we have discussed already in Theorem~\ref{t:main} and its proof,
there is a large number of sporadic cases which have the real cohomology of a product of spheres but different integral cohomology. All these spaces 
are geometrically formal by Proposition~\ref{prod}. To end this section, we discuss the diffeomorphism type of one of these sporadic examples.

Consider the composition of standard inclusions $SO(3)\subset SU(3)\subset SU(4)$. The quotient
$X=SU(4)/SO(3)$ has the real cohomology of $S^5\times S^7$, but different integral cohomology.
\begin{prop}
The manifold $X=SU(4)/SO(3)$ is not diffeomorphic to a non-trivial product or to a symmetric space.
\end{prop}
\begin{proof}
Suppose that $X$ were diffeomorphic to some non-trivial product. Then the factors would have to be simply connected real homology spheres of 
dimensions $5$ and $7$ respectively. As we are assuming that $X$ is diffeomorphic
to the product, not just homotopy equivalent to it, the factors must in fact be homogeneous spaces. Therefore, all the candidates must occur in the 
classification of Onishchik~\cite{On,OR} and Kramer~\cite{Kramer} .
In dimension $5$ the candidates are $S^5$ and $SU(3)/SO(3)$, and in dimension $7$ they are $S^7$, $V_2(\R^5)$ and $Sp(2)/Sp(1)$ with a non-standard
embedding of the subgroup. Looking at the third homotopy groups we have $\pi_3(X)=\pi_3(SU(3)/SO(3))=\Z_4$, but $\pi_3(V_2(\R^5))=\Z_2$ and $\pi_3(Sp(2)/Sp(1))=\Z_{10}$,
cf.~\cite[p.~65]{Kramer}. Thus, on the one hand, the only product of homogeneous spaces that has the same third homotopy group as $X$ is $S^7\times SU(3)/SO(3)$. 
On the other hand, $X$ and $S^7\times SU(3)/SO(3)$ have different sixth homotopy groups.
This is so because the exact homotopy sequence of the fibration $SU(4)\longrightarrow X$ shows that $\pi_6(X)=\pi_5(SO(3))$, which is of order $2$, whereas
$\pi_6(SU(3)/SO(3))$ is known to be of order $4$, see~\cite{L}.

Finally if $X$ were diffeomorphic to a symmetric space, then by what we just proved, that symmetric space would have to be irreducible. However,
by the classification of irreducible symmetric spaces, there is no such space with the real cohomology of $S^5\times S^7$.
\end{proof}
Note that if we discuss only the homotopy type of $X=SU(4)/SO(3)$, then we have to consider products of real homology spheres which are not necessarily 
homogeneous, and the argument breaks down.

\section{Biquotients}\label{s:biq}

Let $G$ be a compact Lie group, and $U$ a closed subgroup of $G\times G$
which acts freely on $G$ by $(u_1,u_2)\cdot g=u_1gu_2^{-1}$.
Then the orbit space is a smooth manifold denoted $G/U$, and is called
a biquotient of $G$. A biinvariant metric on $G$ descends to a metric
on any biquotient. More generally, one can consider subgroups
$U\subset \rm{ISO}(G)$ acting freely. In some cases $\rm{ISO}(G)$ is
strictly larger than $G\times G$, so that one gets more examples.

Biquotients were first studied systematically by Eschenburg~\cite{E0,E} as a
source of examples of manifolds with positive sectional curvature.
Among Eschenburg's biquotients in~\cite{E} there are several examples for
which we can prove easily that they are not geometrically formal.

\begin{ex}
There is a $6$-dimensional example $M=G/U$ whose real cohomology is very
similar but not isomorphic to that of $SU(3)/T^{2}$, considered in~\cite{KT}. Moreover, $M$ and $SU(3)/T^{2}$
have the same integral homology groups and the same cohomology algebra mod $2$.
We can show that this biquotient is not geometrically formal by using almost the
same argument as the one used for $SU(3)/T^{2}$ in~\cite{KT}.

Let $G = U(3)$ and
$$
U\subset\{ (D(a, a,\bar{a}), D(b, c, 1)); a, b, c\in S^{1}\}\subset G\times
G \ .
$$
Eschenburg~\cite{E} proved that $H^{*}(G/U)$ is generated by two
generators $x$ and $y$ of degree $2$ with the relations $xy = y^{2} - x^{2}$ and
$x^{3} = 0$. Setting $z=\frac{1}{\sqrt{5}}(x-2y)$, we get $z^{2}=x^{2}$ and
$z^{3}=-\frac{2}{\sqrt{5}}xy^{2}\neq 0$. If $G/U$ is
geometrically formal, then it has two $2$-forms $x$ and
$z$, such that $x^{3}=0$, $x^{2}=z^{2}$ and $z^{3}$ is volume form.

If $G/U$ were geometrically formal, then for a formal metric these relations would hold 
pointwise for the harmonic forms representing the cohomology classes $x$ and $z$
and their products. The form $x$ would then have rank $4$ everywhere. If $v$ were a vector 
in its kernel, the relation $x^2 = z^2$ would show $i_vz\wedge z=0$, which would contradict the fact
that $z^3$ would have to be a volume form.
\end{ex}

The second example uses a different argument to obstruct geometric formality, related to symplectic
structures defined by harmonic forms of formal metrics.
\begin{ex}
Let $G=SU(3)$ and
$$
U=T^{1}\times T^{1}=\{D(a^{k},a^{l},a^{-k-l}),
D(b^{m},b^{n},b^{-m-n}); a,b\in S^{1} \} \ .
$$
Using the results of~\cite{E} on the cohomology of biquotients, it is easy
to compute $H^{*}(G/U)$. We get the following algebra
structure: there are two linearly independent generators $x$ and $y$
in degree $2$, subject to the relations $x^{2}=y^{2}$ and
$x^{3}=y^{3}$. Then $H^{4}$ is spanned by $xy$ and by $x^{2}=y^{2}$,
and $H^{6}$ is spanned by $x^{3}=y^{3}=x^{2}y=xy^{2}$.
If we assume that $M$ is geometrically formal, then the harmonic form
representing $x+y$ is a symplectic form. It then follows from
$(x-y)(x+y)=0$ that $x-y$ vanishes, which contradicts the linear
independence of $x$ and $y$, because on a symplectic six-manifold
the wedge product with the symplectic form is an isomorphism between
two-forms and four-forms. Thus $M$ cannot be geometrically formal.
\end{ex}

Next we consider Totaro's biquotients of $S^3\times S^3\times S^3$, which he studied in~\cite{tot} as an 
example of a family of $6$-manifolds with nonnegative sectional curvature, but with infinitely many distinct
isomorphism classes of rational cohomology rings. For all these manifolds we
prove that they can not be geometrically formal because of the structure of their cohomology rings.
\begin{thm}
Totaro's biquotients~\cite{tot} are not geometrically formal.
\end{thm}
\begin{proof}
The rational cohomology ring of the biquotients Totaro considers
has three  generators $x_1, x_2, x_3$  in degree $2$ satisfying the  relations
\begin{equation}\label{rel}
\begin{split}
x_1^2 &=0 \\
ax_1x_2+x_2x_3+x_2^2 &=0 \\
bx_1x_3 + 2x_2x_3 +x_3^2 &=0
\end{split}
\end{equation}
for some integers $a$ and $b$, see~\cite{tot}. If we assume that these manifolds are geometrically
formal, then the same relations~\eqref{rel} hold at the level of their representative
harmonic forms, which we denote by the same letters.
In order to prove that geometric formality leads to a contradiction, we differentiate
the following cases according to the values of the constants $a$ and $b$.

{\bf 1.} We first consider the case when $a\cdot b\neq 0$. Then we
can obviously normalize $x_1$ such that $a$ becomes $1$.
Let us consider the forms $y_1=x_1+\frac{3}{b}x_2$ and
$y_2=x_1+\frac{3}{2}x_3$. Using the relations~\eqref{rel} we see
that $y_1^3=y_2^3=0$ and $x_1y_1^2$, $x_1y_2^2$ and $y_1y_2^2$ are
volume forms.

Since the cubes of $y_1$ and $y_2$ vanish, it follows that
$\dim Ker (y_1)=\dim Ker (y_2)=2$, and since $y_1y_2^2$ is a volume
form on $M$, it follows that $Ker (y_1)\cap Ker (y_2)=0$.

If we rewrite the relations~\eqref{rel} in terms of $x_1$, $y_1$ and
$y_2$, we find, after some straightforward calculations, that
\begin{equation*}
\begin{split}
x_1^2 &=0 \\
(1-2b)x_1y_1-2x_1y_2+2y_1y_2+by_1^2 &= 0 \\
-2bx_1y_1+(b-4)x_1y_2+2by_1y_2+2y_2^2 &=0 \ .
\end{split}
\end{equation*}
If we multiply the first relation by $b-4$ and add it to the second one multiplied by
$2$, we get
\begin{equation}\label{nrel}
(5b-2b^2-4)x_1y_1+(6b-8)y_1y_2+b(b-4)y_1^2+4y_2^2=0 \ .
\end{equation}

Let $u_1\in Ker (y_1)$ and $u_2\in Ker (y_2)$ such that
$u_2\notin \R \{u_1 \} \oplus Ker (x_1)$ (such a $u_2$ always exists,
otherwise we would have $Ker (x_1)\cap Ker (y_2)\neq 0$ which is
a contradiction with the fact that $x_1y_2^2$ is a volume form).
In other words $T_{p}M = \{u_1, u_2\}\oplus Ker (x_1)$.
Then if we contract the  equation~\eqref{nrel} with $u_1$ and $u_2$ we get
\begin{equation}\label{nnrel}
(5b-2b^2-4)i_{u_2}(y_1\wedge i_{u_1}x_1)+(6b-8)i_{u_2}y_1\wedge i_{u_1}y_2=0 \ .
\end{equation}
Look at the $1$-forms $i_{u_2}y_1$ and $i_{u_1}y_2$. Since the
dimension of the kernel of each of them is $\geq 5$, it follows that
the intersection $L$ of their kernels has dimension $\geq 4$.
Since $x_1^2=0$, it follows that $\dim Ker (x_1)=4$, and thus
$\dim (Ker (x_1)\cap L) \geq 2$.

Let $w\in Ker (x_1)\cap L$.
If we contract~\eqref{nnrel}  with $w$ we get
$$
x_1(u_1, u_2)i_{w}y_1 = 0,
$$
since the coefficient $5b-2b^2-4$ has no real zeros.

This is in contradiction with the fact that $x_1y_1^2$ is a volume
form. Namely, take $w_1, w_2, w_3 \in Ker (x_1)$ such that together with $u_1,
u_2, w$ they form a basis of $T_{p}M$. Then we have that
$$
x_1y_1^2(u_1,u_2,w,w_1,w_2,w_3)=0 \ ,
$$ 
which is impossible.

{\bf 2.} If $a=0$ and $b\neq 0$, we first normalize $x_1$ to get $b=1$
and proceed  as in the previous case taking $y_1=x_1+3x_2$ and
$y_2=x_1+6x_3$.

{\bf 3.} If $a\neq 0$ and $b=0$, we again first normalize $x_1$ to
have $a=1$ and take $y_1=x_2$ and $y_2=x_1+\frac{3}{2}x_3$.

{\bf 4.} For $a=b=0$ we have the following relations in cohomology
$$
x_1^2=0,\;\; x_2x_3+x_2^2=0,\;\; 2x_2x_3+x_3^2=0 \ .
$$
Take $y_1=x_2+x_3$ and $y_2=x_2+\frac{1}{2}x_3$. The cohomology relations implies that 
$y_1^3=y_2^3=0$  and that $y_1y_2^2$ is a volume form on $M$. If we rewrite the above relations in terms of 
$y_1$ and $y_2$, we obtain the following
$$
y_2^2-y_1y_2=0,\;\; y_1^2-2y_1y_2=0 \ .
$$
Now take $u\in Ker (y_1)$. The relations imply $i_{u}y_2^2=0$, which contradicts the  fact
that $y_1y_2^2$ is a volume form.
\end{proof}

\section{Two-sphere bundles over the complex projective plane}\label{S2-bundles}

For our calculations in~\cite{KT}, the example $SU(3)/T^{2}$ was the
crucial case, from which all others were derived by various
generalisations. The inclusions $T^{2}\hookrightarrow
U(2)\hookrightarrow SU(3)$ show that $SU(3)/T^{2}$ fibers over $\C
P^{2}$ with fiber $S^{2}$. In fact, this fibration is well known as
the twistor fibration of $\C P^{2}$, compare~\cite{KT2}. We now generalize to
arbitrary $2$-sphere fibrations over $\C P^{2}$.

Assume first that we have an arbitrary smooth oriented fibration with
fiber $S^{2}$. Then because the orientation-preserving diffeomorphism group
of $S^{2}$ is homotopy equivalent to $SO(3)$, we may assume that the
structure group is $SO(3)$. Then $M$ is the unit sphere bundle in the
associated rank $3$ vector bundle $V$.

As $SO(3)$ coincides with the
projective unitary group $PU(2)$, every $2$-sphere bundle is the
projectivisation of a complex rank $2$ vector bundle $E$. One can
recover $V$ by passing to the adjoint bundle. We have
\begin{alignat*}{1}
    w_{2}(V) &=c_{1}(E) \ \ \pmod 2 \ , \\
    p_{1}(V) &=c_{1}^{2}(E)-4c_{2}(E) \ .
    \end{alignat*}
In many cases, for example when $w_{2}(V)$ is non-trivial and when $p_{1}(V)$ is divisible by $8$, it follows 
from a result of Grove and Ziller~\cite[Theorem C]{GZ} that $M$ admits a cohomogeneity 
one action of $SO(3)\times SO(3)$. Thus these manifolds are highly symmetric, and are, in
this sense, close to homogeneous.
\begin{thm}
    Let $M^{6}$ be the total space of an $S^{2}$-bundle $E$ over $\C P^{2}$.
    Then $M$ is geometrically formal if and only if it is the trivial bundle $S^{2}\times\C P^{2}$.
    \end{thm}
    \begin{proof}
Assume $M^{6}$ is the total space of the projectivisation of a
complex
rank $2$ vector bundle $E$ over $\C P^{2}$. Then the cohomology of
$M$ is generated multiplicatively by two degree $2$ classes $x$ and
$y$, where we use $x$ for the generator pulled back from $\C P^{2}$
and $y$ for a class which restricts as a generator to every fiber. By
the definition of Chern classes we choose $y$ so that
\begin{equation}\label{Chern}
    y^{2}+c_{1}(E)xy+c_{2}(E)x^{2}=0 \ ,
    \end{equation}
where by an obvious abuse of notation we use $c_{i}(E)$ for the Chern
numbers $\langle c_{i}(E),[\C P^{i}]\rangle$. We also have $x^{3}=0$.

After replacing $y$ by a linear combination of $x$ and $y$ we may
assume that
\begin{equation}\label{Chern2}
    y^{2}+c x^{2}=0 \ ,
    \end{equation}
    where $c=-\frac{1}{4}p_{1}(V)$ vanishes if and only of $M$ is the
    trivial bundle. To pass from~\eqref{Chern} to~\eqref{Chern2} we
    can twist $E$ by a line bundle.
    This does not affect the projectivisation, but we can kill the
    first Chern class if $w_{2}(V)=0$. If this is not the case, we
    can
    still make the required base change, which amounts to twisting by
    a virtual line bundle whose Chern class is half-integral. In this
    case the constant $c$ is not integral.

    Now~\eqref{Chern2} together with $x^{3}=0$ implies $xy^{2}=0$. We
    also have $y^{3}+cx^{2}y=0$. If $c\neq 0$, we conclude that
    $y^{3}\neq 0$, for otherwise there would be no degree $6$
    cohomology in $M$. Thus, if $c\neq 0$ and $M$ is geometrically
    formal, then the harmonic $2$-form representing $y$ is
    nondegenerate. But the harmonic form representing $x$ has a
    nontrivial kernel as $x^{3}=0$, and if we contract~\eqref{Chern2}
    with an element $v$ in this kernel, we deduce $2i_{v}y\wedge
    y=0$, which contradicts the nondegeneracy of $y$.

    Thus we have proved that a nontrivial $M$ cannot be geometrically
    formal. Conversely, the trivial bundle $S^{2}\times\C P^{2}$ is
    geometrically formal with respect to a product of K\"ahler
    metrics, cf.~\cite{K}. In fact, it is a symmetric space.
\end{proof}

\bigskip

\bibliographystyle{amsplain}

\begin{thebibliography}{10}

\bibitem{Adams}
J.~F.~Adams, {\em On the non-existence of elements of Hopf invariant one},
Ann.~of Math.~{\bf 72} (1960), 20--104.

\bibitem{AW}
S.~Aloff and N.~R.~Wallach, {\em An infinite family of distinct
$7$-manifolds admitting positively curved Riemannian structures},
Bull.~Amer.~Math.~Soc.~{\bf 81} (1975), 93--97.

\bibitem{Besse}
A.~L.~Besse, {\sl Einstein Manifolds}, Springer Verlag 1987.

\bibitem{Borel}
A.~Borel, {\em Les bouts des espaces homog\`enes de groupes de Lie}, Ann.~Math.~{\bf
58} (1953), 443--457.

\bibitem{Borel3}
A.~Borel, {\em La cohomologie$\mod 2$ de certains espaces  homog\`enes}, 
Comment.~Math.~Helv.~{\bf  27} (1953), 165--197.

\bibitem{BH}
A.~Borel and F.~Hirzebruch, {\em Characteristic classes and
homogeneous spaces, I}, Amer.~J.~Math.~{\bf 80} (1958), 459--538.

\bibitem{BS}
A.~Borel and J.-P.~Serre, {\em Groupes de Lie et puissances r\'eduites de Steenrod}, Am.~J.~of 
Math.~{\bf 75} (1953), 409--448.

\bibitem{Bott}
R.~Bott, {\sl Lectures on $K(X)$}, W.~A.~Benjamin Inc., New York, Amsterdam 1969.

\bibitem{DGMS}
P.~Deligne, P.~Griffiths, J.~W~Morgan and D.~Sullivan, {\em Real homotopy theory of K\"ahler manifolds},
Invent.~Math.~{\bf 29} (1975), 245--274. 

\bibitem{E0}
J.~H.~Eschenburg, {\em Inhomogeneous spaces of positive curvature}, Differential Geom.~Appl.~{\bf 2} (1992), 123--132.

\bibitem{E}
J.~H.~Eschenburg, {\em Cohomology of biquotients}, Manuscripta
Math.~{\bf 75} (1992), 151--166.

\bibitem{G}
M.~E.~Gilmore, {\em Complex Stiefel manifolds, some homotopy groups and vector fields}, Bull.~Amer.~Math.~Soc.~{\bf 73} (1967), 630--633.

\bibitem{GZ}
K.~Grove and W.~Ziller, {\em Lifting group actions and nonnegative curvature}, Preprint arXiv:0801.0767 v1 [math.DG] 5 Jan 2008.

\bibitem{HS}
W.-Y.~Hsiang and J.~C.~Su, {\em On the classification of transitive effective actions on Stiefel manifolds},
Trans.~Amer.~Math.~Soc.~{\bf 130} (1968), 322--336.

\bibitem{Huse}
D.~Husemoller, {\sl Fibre Bundles}, Second Edition, Springer Graduate Text in Mathematics {\bf 20}, Springer Verlag, New York, Heidelberg, Berlin 1975.

\bibitem{James}
I.~M.~James, {\em On the homotopy type of Stiefel manifolds}, Proc.~Amer.~Math.~Soc.~{\bf 29} (1971), 151--158.

\bibitem{JW}
I.~M.~James and J.~H.~C.~Whitehead, {\em The homotopy theory of sphere bundles over spheres (I)}, 
Proc.~London Math.~Soc.~{\bf 4} (1954), 196--218.
 
\bibitem{KN}
S.~Kobayashi and K.~Nomizu, {\sl Foundations of Differential Geometry},
Interscience Publishers 1963.

\bibitem{K}
D.~Kotschick, {\em On products of harmonic forms}, Duke
Math.~J.~{\bf 107} (2001), 521--531.

\bibitem{KT}
D.~Kotschick and S.~Terzi\'c, {\em On formality of generalised symmetric spaces},
Math.~Proc.~Cam.~Phil.~Soc.~{\bf 134} (2003), 491--505.

\bibitem{KT2}
D.~Kotschick and S.~Terzi\'c, {\em Chern numbers and the geometry of partial flag manifolds},
Comment.~Math.~Helv. (to appear).

\bibitem{Kramer}
L.~Kramer, {\sl Homogeneous Spaces, Tits Buildings, and Isoparametric Hypersurfaces}, Memoirs of the Amer.~Math.~Soc.~Vol.~{\bf 158} Number 752 (2002).

\bibitem{Doan Kuin'}
D.~Kuin', {\em The Poincar\'e polynomials of compact homogeneous
Riemannian spaces with irreducible stationary group} (Russian), Trudy
Sem.~Vektor.~Tenzor.~Anal.~{\bf 14} (1968), 33--93.

\bibitem{L}
A.~T.~Lundell, {\em Concise tables of James numbers and some homotopy of classical Lie groups and associated homogeneous spaces}, in
{\sl Algebraic Topology -- Homotopy and Group Cohomology}, ed.~ J.~Aguad\'e, M.~Castellet and F.~R.~Cohen, Springer Lecture Notes in Mathematics
{\bf 1509}, Springer Verlag 1992.

\bibitem{M}
M.~Mimura, {\em The homotopy groups of Lie groups of low rank}, J.~Math.~Kyoto Univ.~{\bf 6} (1967), 131--176.

\bibitem{O}
K.~Oguchi, {\em Homotopy groups of $Sp(n)/Sp(n-2)$}, J.~Fac.~Sci.~Univ.~Tokyo {\bf 16} (1969), 179--201. 

\bibitem{On}
A.~L.~Onishchik, {\em Transitive compact transformation groups}, Amer.~Math.~Soc.~Transl.~{\bf 55} (1966), 153--194.

\bibitem{OR}
A.~L.~Onishchik, {\em Topology of transitive transformation groups} (Russian), Fizmatlit Nauka Moscow 1995.

\bibitem{Sem}
U.~Semmelmann, Private communication.

\bibitem{tokyo}
D.~Sullivan, {\em Differential Forms and the Topology of Manifolds},
in {\sl Manifolds Tokyo 1973}, ed.~A.~Hattori, University of Tokyo
Press 1975.

\bibitem{Takeuchi}
M.~Takeuchi, {\em On Pontrjagin classes of compact symmetric spaces},
J.~Fac.~Sci.~Univ.~Tokyo Sec I {\bf 9} (1962), 313--328.

\bibitem{tot}
B.~Totaro, {\em Curvature, diameter, and quotient manifolds},
Mathematical Research Letters {\bf 10} (2003), 191--203.

\bibitem{T}
S.~Terzi\'c, {\em Cohomology with real coefficients of generalized
symmetric spaces} (Russian), Fundam.~Prikl.~Mat.~Vol. {\bf 7} no.1
(2001), 131--157.

\bibitem{T1}
S.~Terzi\'c, {\em Rational homotopy groups of generalised symmetric
spaces}, Math.~Z.~{\bf 243} (2003), 491-523.

\bibitem{WZ}
M.~Wang and W.~Ziller, {\em Symmetric spaces and strongly isotropy
irreducible spaces}, Math.~Ann.~{\bf 296} (1993), 285--326.

\bibitem{W-G}
J.~A.~Wolf and A.~Gray, {\em Homogeneous spaces defined by Lie group
automorphisms. I} and {\em II}, J.~Differential Geometry {\bf 2} (1968), 77--114
and 115--159.

\end{thebibliography}

\bigskip

\end{document}